\documentclass[reqno]{amsart}

\usepackage[all]{xy}
\usepackage{enumerate}
\usepackage{amscd}
\usepackage{amssymb,amsmath}
\usepackage{graphicx}

\theoremstyle{plain}
  \newtheorem{theorem}{Theorem}[section]
  
  \newtheorem*{corollary*}{Corollary}
  
  \newtheorem*{lemma*}{Lemma}
  
\theoremstyle{definition}
  \newtheorem{definition}[theorem]{Definition}
  
  \newtheorem{observation}[theorem]{Observation}

\usepackage{amssymb,amsmath}

\begin{document}

\title{ Topology of real  and angle valued maps \\
 and graph representations. 
}


\author{Dan Burghelea}

\address{Dept. of Mathematics,
         The Ohio State University,
         231 West 18th Avenue,
         Columbus, OH 43210, USA.}

\email{burghele@mps.ohio-state.edu}

\thanks{. 
  The  author acknowledge partial support from NSF grant MCS 0915996.
  }



\date{\today}

\begin{abstract}

Using graph representations a new class of computable topological invariants
associated with a tame real or angle valued map were recently introduced, 
providing a theory which can be viewed as an alternative to Morse-Novicov theory for
real or angle valued Morse maps.
The invariants are "barcodes" and "Jordan cells".
From them one can derive all familiar topological invariants which can be derived via
Morse-Novikov theory, like  the Betti numbers  and in the case of angle  valued maps  also the Novikov Betti numbers and the monodromy.  
Stability results for bar codes and the homotopy invariance of
the Jordan cells  are the key results, and two new polynomials  for any $r$  associated to
a continuous nonzero complex valued map provide  potentially interesting refinements of the Betti numbers and  of the Novikov Betti numbers.
In our theory the bar codes which are  intervals with ends critical values/angles, the Jordan cells and the " canonical long exact sequence" of a tame map are the analogues
of instantons between rest points, closed trajectories and of the Morse-Smale complex of the gradient of a
Morse function in the Morse-Novikov theory.
\end{abstract}

\maketitle

\setcounter{tocdepth}{1}

\vskip .2in 

\section  {Introduction}.
\vskip .2in

This is essentially the lecture delivered at the Congress of the Romanian mathematicians, Brasov, June 2011 under the title "New topological invariants for angle valued maps".

The presentation summarizes work  done in \cite{BD11}, \cite {BH12} and \cite {B12}. 
Using graph representations  and inspired by persistence theory  \cite {ELZ02} \cite {CSD09}  and \cite {BD11}
a new class of computable topological invariants
associated with a tame real or angle valued maps were recently introduced, providing 
a theory which can be viewed as an alternative to Morse-Novicov theory for
real or angle valued Morse maps.
The invariants are "barcodes" and "Jordan cells" and,  when
the underlying space is a simplicial complex and the map is simplicial, can be calculated
by
algorithms of the same complexity as the ones which calculate the Betti numbers.
From them one can derive all familiar topological invariants which can be derived via
Morse-Novikov theory and a few more. Stability results for bar codes, Theorem \ref{T3}, and homotopy invariance of
the Jordan cells  and of the cardinality of some sets of bar codes, Theorems \ref{T1}, \ref{T2}, are the key results, and two new polynomials  associated to
a continuous nonzero complex value map provide  potentially interesting refinements of 
the Betti numbers and  of the Novikov-Betti numbers.

In our theory the bar codes which are intervals with ends critical values/angles, the 
Jordan cells and the canonical long exact sequence of a tame map are the analogues
of the instantons, the closed trajectories and the Morse Smale complex of the gradient of a
Morse function in the Morse Novikov theory.

Note that  almost all maps are tame in particular all Morse maps on a smooth manifold or on a stratified space and all simplicial maps on a simplicial complex. In the case of angle (circle)  valued maps the space of tame maps have the same homotopy type as the space of all 
continuous maps.  This  is not the case of Morse maps. Note also that in case of Morse angle  valued maps the cardinality of instantons and closed trajectories might not
be finite, but  the cardinality of the set of bar codes and Jordan cells is always finite.  This presentation contains  Theorem \ref{T3} which was not present in my lecture at the Romanian congress. The organization of the material  is  also slightly different.

\vskip .2in
\section  {Topology}
\vskip .2in
Let $\kappa$ be a field and $\overline \kappa$ its algebraic closure. 
Denote by $\kappa[t, t^{-1}]$ the ring of the Laurent  polynomials and by  $\kappa[[t, t^{-1}]$
the field of Laurent power series with coefficients in  $\kappa.$  Clearly $\kappa[t, t^{-1}]\subset \kappa[[t, t^{-1}].$

Let $X$ be a compact ANR. Denote by $H_r(X)$ the singular homology with coefficients in $\kappa$ and call
$$\beta_r(X)= \dim H_r(X)$$ the $r-$th 
 {\it Betti number}  of $X.$ 

Let $\xi\in H^1(X,\mathbb Z).$ For the pair $(X,\xi)$ consider $\tilde X\to X$ the infinite cyclic cover associated with $\xi,$ precisely the pull back of the canonical infinite cyclic cover $\mathbb R\to S^1= \mathbb R/\mathbb Z$ by a map $f: X\to S^1$ representing $\xi.$ Let $T:\tilde X\to \tilde X$  be the deck transformation.

Note that $H_r(\tilde X)$ is a $\kappa [t, t^{-1}]$ module with the multiplication by $t$ induced by the deck transformation $T.$  Consider $NH_r(X,\xi)= H_r(\tilde X)\otimes_{\kappa [t, t^{-1}]} \kappa[[t, t^{-1}].$  This is the Novikov homology which is a vector space over the field $\kappa [[t, t^{-1}].$ We call 
$$\beta N_r(X;\xi)= \dim_{\kappa [[t, t^{-1}]} NH_r(X;\xi)$$ the $r-$th {\it Novikov-Betti number} of $(X,\xi)$.

Let $$ V(\xi):=\ker \{H_r(\tilde X)\to NH_r(X;\xi)\}$$ be the kernel of the linear map induced by tensoring $H_r(\tilde X)$ with  $\kappa[t, t^{-1}]\to \kappa[[t, t^{-1}]$ over the ring $\kappa[t, t^{-1}].$
The $\kappa[t, t^{-1}]$ module $V(\xi)$ is a finite dimensional vector space over the field $\kappa.$ The multiplication by $t$ can be viewed as a $\kappa-$ linear isomorphism   $T(\xi):V(\xi)\to V(\xi).$ 
The pair $(V(\xi),T(\xi))$ will be referred to as the {\it monodromy}  associated with $\xi.$

\vskip .2 in

\section  {Tame maps}.
\vskip .1in
 \begin{definition}
A continuous map $f\colon X\to\mathbb R$ resp.\ $f\colon X\to S^1$, $X$ a compact ANR, is \emph{tame} if the following hold:
\begin{enumerate}
\item
Any fiber $X_\theta=f^{-1}(\theta)$ is the deformation retract of an open neighborhood. 
\item
Away from a finite set of numbers resp. angles $\Sigma=\{s_1, s_2, \cdots s_r\} \subset \mathbb R$, resp.\ $S^1$
the restriction of $f$ to $X\setminus f^{-1}(\Sigma)$ is a fibration.
\end{enumerate}
\end{definition}

For any  real  resp. angle valued tame map we have the finite set of numbers  $ s_1< s_2 <\cdots s_{N-1} <s_N $   resp. $ 0\leq \theta_1< \theta_2 <\cdots \theta_{m-1} <\theta_m <2\pi$  where the homotopy type of the fibers change. The numbers $s_i$ resp. $\theta_i$  are  the {\it critical values} of the tame map $f.$

In the case of a real valued map the Betti numbers and in the case of an angle valued map  the Betti numbers, the Novikov-Betti numbers and the monodromy  
can be recovered from the  invariants associated with the tame map. These invariants  are the {\it bar codes } and the {\it Jordan cells} and are computable, cf section 10. 
\vskip .2in

\section {Bar codes and Jordan cells. The invariants of a tame map $f.$} 
 \vskip .1in
 Bar codes are finite intervals $I$ of real numbers of four types: \begin {enumerate}
 \item Type 1, closed , $[a,b]$ 
  with $a\leq b, $ 
\item Type 2, open, $(a,b)$ 
  with $a< b, $ 
\item Type 3, left open right closed, $(a,b]$ 
  with $a< b, $ 
  \item Type 4, left closed right open,  $[a,b)$ 
  with $a<b. $ 
\end{enumerate}

Jordan cells are pairs $J = (\lambda \in \overline \kappa, n\in \mathbb Z_{\geq 0}).$
A Jordan cell  should be interpreted as a 
 a matrix 

\begin{equation}\label{E1}
T(\lambda,k)=
\begin{pmatrix}
\lambda & 1       & 0      & \cdots  & 0      \\
0       & \lambda & 1      & \ddots  & \vdots \\
0       & 0       & \ddots & \ddots  & 0      \\
\vdots  & \ddots  & \ddots & \lambda & 1      \\
0       & \cdots  & 0      & 0       & \lambda
\end{pmatrix}_.
\end{equation}



\vskip .1in

For a tame map $f:X\to \mathbb R$   and any integer $r$ we associate cf. section 6 a collection of
of bar codes $$\mathcal B_r(f)= \mathcal B^c_r (f) \sqcup  \mathcal B^o_r(f) \sqcup  \mathcal B^{o,c}_r(f) \sqcup  \mathcal B^{c,o}_r(f)$$ 
with $\mathcal B^c, \mathcal B^o,\mathcal B^{o,c},\mathcal B^{c,o}$ of type 1,2,3,4. whose ends $a,b$ are  critical values. 

For a tame map $f:X\to \mathbb S^1$   and any integer $r$ we associate a collection of
of bar codes 
$$\mathcal B_r(f)= \mathcal B^c_r (f) \sqcup  \mathcal B^o_r(f) \sqcup  \mathcal B^{o,c}_r(f) \sqcup  \mathcal B^{c,o}_r(f)$$ 
with $\mathcal B^c, \mathcal B^o,\mathcal B^{o,c},\mathcal B^{c,o}$ of type 1,2,3,4  with ends $a,b$ and
Jordan cells 
$\mathcal J(r),$ $J=(\lambda(J), k(J)),$ $\lambda(J)\in \overline \kappa,$ $k(J)\in \mathbb Z_{\geq 0}.$
The ends $a,b$ are the first a critical angle $\theta _i$ the second of the form $\theta_j +2\pi k$, $\theta_j$ a critical angle $k$ a non negative integer. 

In both cases (real and angle valued maps)  it will be convenient to record the bar codes  $\mathcal B^c_r(f)\sqcup \mathcal B^o_{r-1}(f) $ as a configuration $C_r(f)$ of points   
in the plane $\mathbb R^2$ resp. the cylinder $\mathbb T$ defined by 
$\mathbb T=\mathbb R^2/\mathbb Z$ or equivalently  $\mathbb C\setminus 0,$ see picture below.   
Precisely $\mathbb T$ is the quotient space of $\mathbb R^2$ the Euclidean plane, by the additive group of integers $\mathbb Z$, w.r.\ 
to the action $\mu: \mathbb Z\times \mathbb R^2 \to \mathbb R^2$ given by $\mu(n;(x,y))=(x+2\pi n,y+2\pi n).$  The identification of $\mathbb T$ to $\mathbb C\setminus 0$ is done via the map $(x,y)\to e^{(x-y) +ix}.$

One denotes by $\Delta \subset \mathbb R^2$ resp. $\Delta \subset \mathbb T$ the diagonal of $\mathbb R^2$ resp. the quotient of the diagonal of $\mathbb R^2$ by the group $\mathbb Z$.
The points above or on diagonal, $(x,y)$, $x\leq y$,   will be used to record closed bar codes $[x,y]\in \mathcal B^c_r(f)$ and 
the points below the diagonal, $(x,y)$, $x>y$ 
to record open bar codes $(y,x)\in \mathcal B^o_{r-1}$.  When $\mathbb T$ is identified to $\mathbb C\setminus 0$ the diagonal $\Delta$ corresponds to the unit circle, the points above or on the diagonal to the the points outside or on the unit circle and those below diagonal to points inside unit circle.

 If we identify a point in $(x,y)\in R^2$ with $z=x+iy,$ hence $\mathbb R^2$ to $\mathbb C,$ and $\mathbb T$ to $\mathbb C\setminus 0$ it is convenient to regard $C_r(f)$ as 
the monic polynomial $P^f_r(z)$ 
whose roots are the elements of $C_r(f)$.  
In the second case $P^f_r(z)$ is 
a monic  polynomial with nonzero free coefficient since  the roots are all nonzero. 
\vskip .2in

\begin{figure}[h]
\includegraphics [height=6cm]{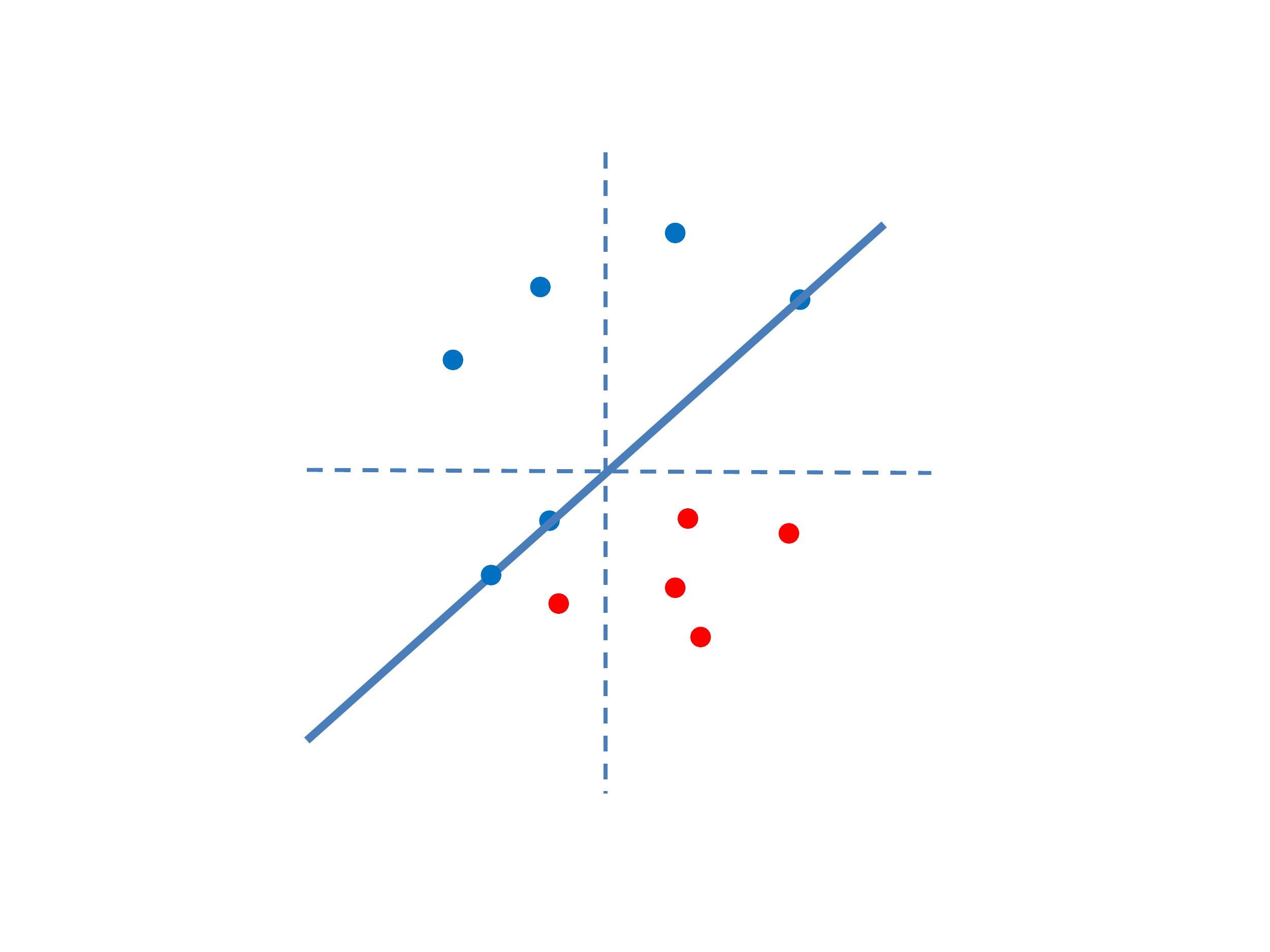}
\includegraphics  [height=6cm] {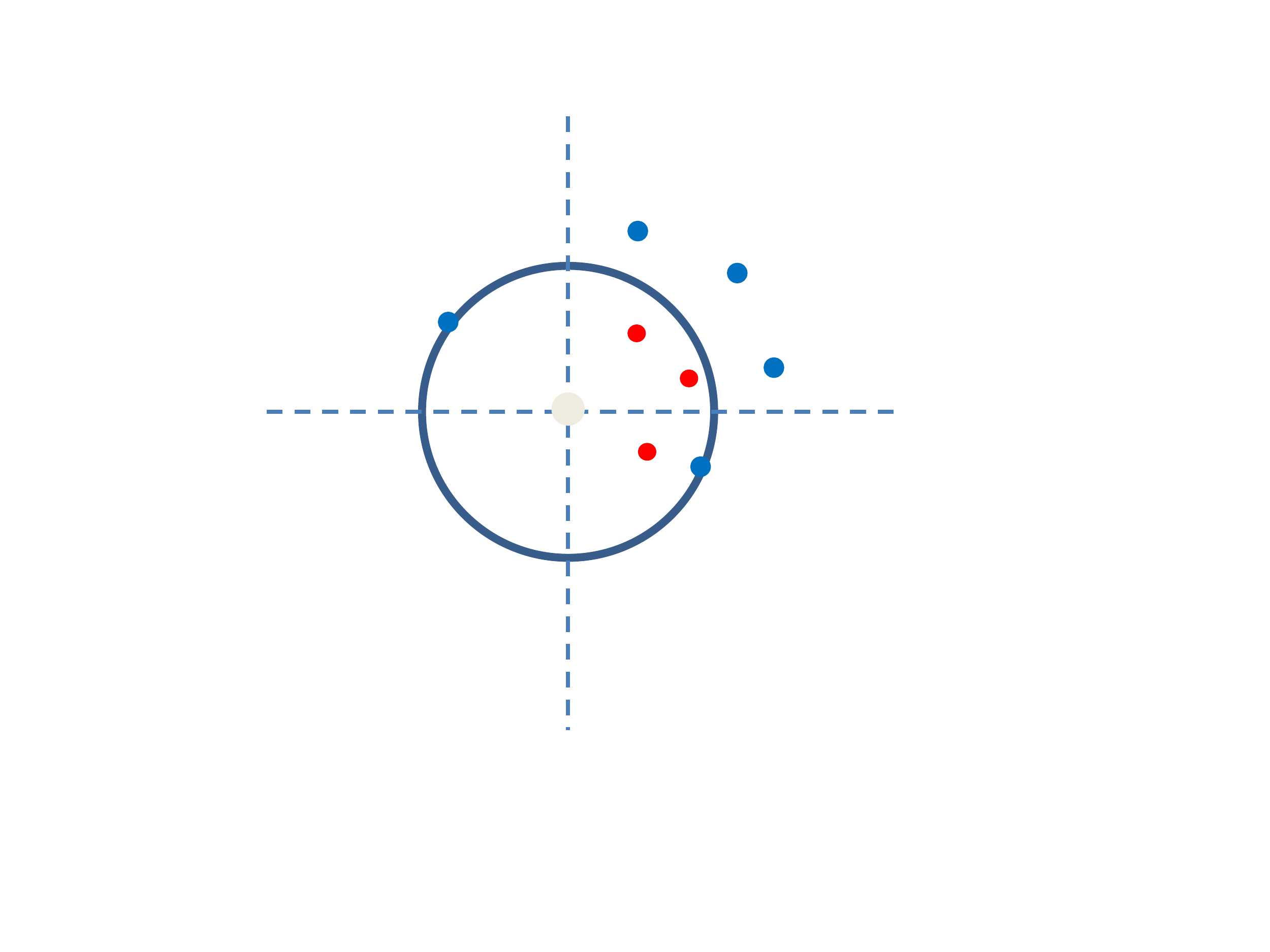}
\caption{Configurations}
\end{figure}

The first  picture is the configuration $Cr_r(f)$ for a tame real valued map  the second for an angle valued map. The points in blue (above or on the diagonal resp. outside or on the unit circle) represent closed $r-$ bar codes, the ones in red (below diagonal resp. inside the unit circle) open $(r-1)-$ bar codes.

\vskip .2in

\section  {Graph representation }
\vskip .1in
To describe  $\mathcal B_r(f)$ and $\mathcal J_r(f)$ we will use two graphs  $\mathcal Z$ for real valued maps, and $G_{2m}$ for angle valued maps. 
The graph $\mathcal Z$ has vertices $x_i$, $i\in\mathbb Z$, and edges $a_i$ from $x_{2i-1}$ to $x_{2i}$ and $b_i$ from $x_{2i+1}$ to $x_{2i}$,

\begin{center}
$$
\xymatrix{
\cdots & x_{2i-1}\ar[l]_-{b_{i-1}} \ar[r]^-{a_i} & x_{2i} & x_{2i+1}\ar[l]_-{b_{i}} \ar[r]^-{a_{i+1}} & x_{2i+2} & \cdots \ar[l]_-{b_{i+1}}
}
$$
The graph $\mathcal Z$
\end{center}
\vskip.1in
\noindent and $\Gamma=G_{2m}$ has vertices  $x_1,x_2,\dotsc,x_{2m}$ and edges $a_i$, $1\leq i\leq m$, and $b_i$, $1\leq i\leq m-1$, as above and $b_m\colon x_1\to x_{2m}$.  
\begin{center}
$$
\xymatrix{ &x_2\\
x_3\ar[ur]_{b_1}\ar[d]^{a_2}& & x_1\ar[ul]^{a_1}\ar[d]_{b_m}\\
x_4& & x_{2m}\\
x_{2m-3} \ar@{<.>}[u]\ar[dr]^{a_{m-1}}& & x_{2m-1}\ar[u]^{a_m}\ar[dl]_{b_{m-1}}\\
& x_{2m-2}
}
$$
The graph $G_{2m}$
\end{center}
.
\vskip .1in
Let $\kappa$ be a fixed field. 
A $\Gamma$-representation $\rho$ is an assignment which to each vertex $x$ of $\Gamma$ assigns a finite dimensional vector space $V_x$  
and to each oriented  arrow  from the vertex $x$ to the vertex $y$ a linear map $V_x\to V_y$. 
The concepts of morphism, isomorphism= equivalence, sum, direct summand, zero and nontrivial representations are obvious.

A $\mathcal Z$-representation is given by the collection 
\begin{equation*}
\rho:=
\begin {cases}
\begin {aligned}  
V_r,\quad \alpha_i:V_{2i-1}\to &V_{2i}, \quad \beta_i:V_{2i+1}\to V_{2i}\\ 
r, i \in &\mathbb Z  \ \ ,
\end{aligned}
\end{cases}
\end{equation*}
while a $G_{2m}$ representation by the collection 
\begin{equation*}
\rho:=
\begin {cases}
\begin {aligned}  
V_r,\quad \alpha_i:V_{2i-1}\to V_{2i}, \quad \beta_i:V_{2i+1}\to V_{2i} \\
1\leq r \leq 2m, \quad 1\leq i\leq m, \quad V_{2m+1}= V_1.
\end{aligned}
\end{cases}
\end{equation*}
Both will be abbreviated by $\rho=\{V_r,\alpha_i,\beta_i \}.$

A finitely supported $\mathcal Z$-representation\footnote{ i.e. all but finitely many vector spaces $V_x$ have dimension zero}, resp.\ an arbitrary $G_{2m}$-representation can be uniquely 
decomposed as a sum of indecomposable representations. In the case of the graph $\mathcal Z$ the indecomposable representations  are indexed by  
one of the four types of intervals (bar codes)  with ends $i,j\in \mathbb Z$, $i\leq j$ for type (1) and $i<j$ for types (2), (3) and (4). For reasons which will be understandable later on we regard the ends $i,j$ 
as associated to the vertices  $x_{2i}, x_{2j}.$
We refer to both the indecomposable representation  and the interval  as \emph{bar code.} 

Here is the description of all {\it bar codes} (for the graph $\mathcal Z$). 

\begin{enumerate}
\item 
$ \rho([i,j]) , i\leq j$  has 
$V_r= \kappa$ for $r = \{2i, 2i+1, \cdots 2j\}$  and $V_r=0$ if $r\ne [2i,2j],$
\item
$ \rho([i,j)) , i< j $ has
$V_r= \kappa$ for $r = \{2i, 2i+1, \cdots 2j\}$  and $V_r=0$ if $r\ne [2i,2j-1],$
\item 
$ \rho((i,j]) , i< j $ has
$V_r= \kappa$ for $r = \{2i, 2i+1, \cdots 2j\}$   and $V_r=0$ if $r\ne [2i+1,2j],$
\item 
$ \rho((i,)]) , i< j $ has
$V_r= \kappa$ for $r = \{2i, 2i+1, \cdots 2j\}$  and $V_r=0$ if $r\ne [2i+1,2j-1],$
\end{enumerate}
 with all $\alpha_i$ and $\beta_i$ the identity  provided that the source and the target are  both non zero.
The above description is implicit in \cite{G72}.

In the case of the graph $G_{2m}$ for simplicity we consider the field $\kappa$ algebraically closed. The indecomposable representations  are indexed by similar intervals (bar codes) with ends $i,j +mk, $ $1\leq i,j\leq m, k\in \mathbb Z_{\leq 0}$, $i\leq j$ with $1\leq i\leq m$ 
and by Jordan cells. Again $i,j$ are associated to vertices $x_{2i}, x_{2j}.$
We refer to both the indecomposable representation 
and the interval  resp. 
the Jordan cell  as \emph{ bar code} resp. 
 \emph {Jordan cell}.

{\bf Type I:} ({\it bar codes} for the graph $G_{2m}$)
For any triple of integers $\{i,j,k\}$, $1 \leq i,j \leq m,$ $k\geq 0$, we have the representations  denoted by 

\begin{enumerate}
\item 
$\rho^I([i, j];k) \equiv \rho^I([i, j+mk]),\quad 1\leq i, j\leq m, k\geq 0$
\item  $\rho^I((i,j];k) \equiv  \rho^I((i, j+mk]),\quad 1\leq i, j\leq m, k\geq 0$ 
\item  $\rho^I([i,j);k) \equiv  \rho^I([i, j+mk)),\quad 1\leq i, j\leq m, k\geq 0$ 
\item $\rho^I((i,j);k)  \equiv  \rho^I((i, j+mk)),\quad 1\leq i, j\leq m, k\geq 0$ 
\end{enumerate}
 described as follows.

Suppose the vertices  of $G_{2m}$ are located counter-clockwise on the unit circle
with  evenly indexed vertices $\{x_2, x_4, \cdots x_{2m}\}$ 
 corresponding  to the angles 
$0< s_1< s_2 <\cdots <s_{m} \leq2\pi.$ 
Draw the spiral curve 
 for $a= s_i$ and $b=s_j+2\pi k$ with the ends a black or an empty circle if the end of the bar code is closed or open 
(see picture below for $k=2$).

\begin{figure}[h]
\includegraphics [height=4.5cm]{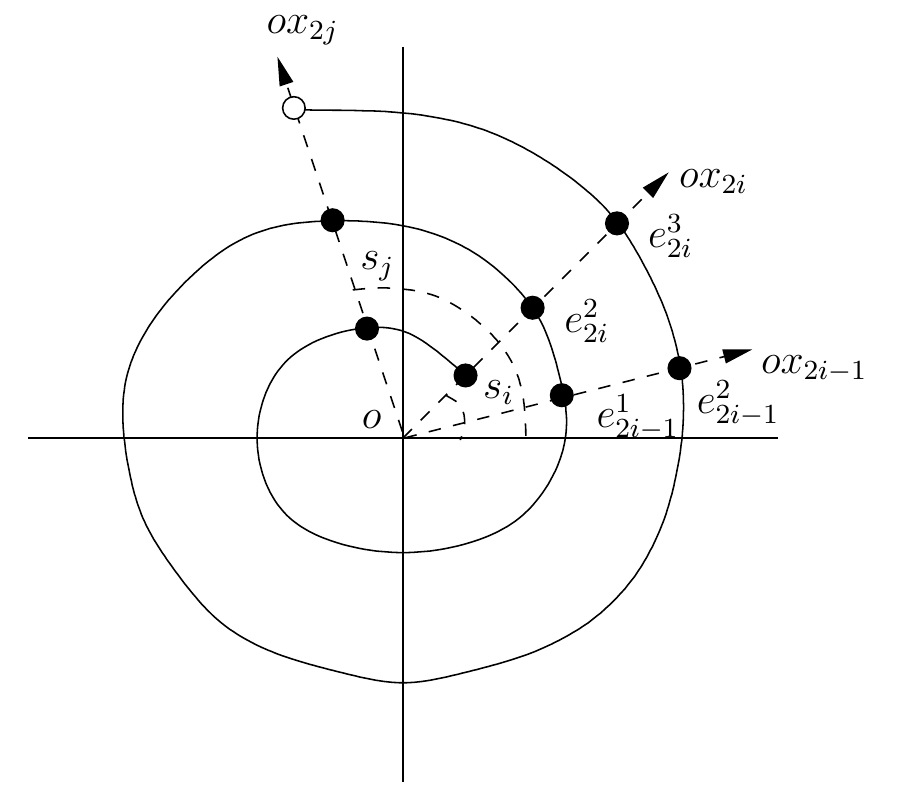}
\caption{The spiral for $[i, j+2m)$.}
\label{barcode}
\end{figure}

Denote by $V_i$ the vector space generated by the intersection points of the spiral with the radius corresponding to the vertex $x_i.$ Let $\alpha_i$ resp.$\beta_i$ be defined as follows:  a generator $e$ of $V_{2i\pm1}$ is sent to the generator $e'$ of $V_{2i}$ if connected by a piece of spiral  and to $0$ otherwise.
\vskip .2in

{\bf Type II: } 
 $\rho^{II}(\lambda, k)$  defined by 

\begin{equation}\label{R1}
\rho^{II}(\lambda, k) = \{V'_r=\lambda^k, \alpha'_1= T(\lambda,k), \alpha'_i=Id\   i\ne 1, \  \beta'_i= Id\}.
\end{equation}

For a $\mathcal Z$-representation or a $G_{2m}$-representation $\rho$ one denotes by $\mathcal B(\rho)$ the set of all 
bar codes and write $\mathcal B(\rho)$ as $\mathcal B(\rho)=\mathcal B^c(\rho)\sqcup \mathcal B^o(\rho)\sqcup \mathcal B^{o,c}(\rho)\sqcup \mathcal B^{c,o}(\rho)$ where 
$\mathcal B^c(\rho)$, $\mathcal B^o(\rho),$ $B^{o,c}(\rho)$ and $ \mathcal B^{c,o}(\rho)$ are the subsets of closed, open, left open right closed, and right open left closed  bar codes.

For a $G_{2m}$ representation $\rho$ one denotes by $\mathcal J(\rho)$  the set of all 
Jordan cells  resp. Jordan cells.

\vskip .2in
\section{The invariants  associated to a tame map $f.$}\label{S6}

Given a tame map $f\colon X\to\mathbb R$ resp.\  $f\colon X\to S^1$ consider the critical values resp.\ the critical angles $\theta_1<\theta_2<\cdots<\theta_m$\footnote {for brevity in writing  denote the critical values of both, real and anglel valued maps by  $\theta_i$.} .  
In the second case  we have $0<\theta_1<\cdots<\theta_m\leq 2\pi$. 
Choose $t_i$, $i=1,2,\dotsc,m$, with $\theta_1<t_1 <\theta_2<\cdots<t_{m-1}<\theta_m<t_m$.  In the second case choose $t_m$ s.t.\ $2\pi<t_m <\theta_1 +2\pi$. 

The tameness of of the map  when $f$ is a real valued map induces the diagram 
\begin{center}
$$
\xymatrix{
\cdots & X_{t_{i-1}}\ar[l]_-{b_{i-1}} \ar[r]^-{a_i} & X_{\theta_i} &X_{t_i}\ar[l]_-{b_{i}} \ar[r]^-{a_{i+1}} & X_{\theta_{i+1}} & \cdots \ar[l]_-{b_{i+1}}
}
$$
\end{center}
\vskip .1in
\noindent and when  $f$ is angle valued map  the diagram
{\scriptsize
\begin{center}
$$
\xymatrix{ &X_{\theta_1}\\
X_{t_1}\ar[ur]_{b_1}\ar[d]^{a_2}& & X_{t_m}\ar[ul]^{a_1}\ar[d]_{b_m}\\
X_{\theta_2}& & X_{\theta_m}\\
X_{t_{m-2}} \ar@{<.>}[u]\ar[dr]^{a_{m-1}}& & X_{t_{m-1}}\ar[u]^{a_m}\ar[dl]_{b_{m-1}}\\
&X_{\theta_{m-1}}
}
_.$$
\end{center}
}
\vskip .1in

Here $X_t= f^{-1}(t)$ resp. $X_\theta= f^{-1}(\theta).$
Different choices of $t_i$ lead to different diagrams but all homotopy equivalent.

For any $r\leq \dim X$ let $\rho_r= \rho(f)$ be the $\mathcal Z$- resp.\ $G_{2m}$-representation associated to the tame map $f$  defined by
$$
V_{2i}= H_r(X_{\theta_i}),V_{2i+1}= H_r(X_{t_{i}}),\quad \alpha_i:V_{2i-1}\to V_{2i},\quad \beta_i:V_{2i+1}\to V_{2i}
$$
with $\alpha_i$ and $\beta_i$ the linear maps induced by the continuous maps $a_i$ and $b_i$ in the diagrams above. 
Here and below $H_r(Y)$ denotes the singular homology in dimension $r$ with coefficients in a fixed chosen field $\kappa .$ 

In order to relate the indecomposable components of $\rho_r$ to the critical values  of $f,$ for a real valued map one converts the intervals $\{i,j\}$  into $\{\theta_i,\theta_j\}$   and for an angle value map the intervals 
$\{i,j+km\}$, $1\leq i,j\leq m$,  into the intervals $\{\theta_i,\theta_j+2\pi k\}$\footnote {we use the symbol "$\{$" for both "$($" and "$[ $" and  "$\}$" for both $")"$ and $" ] ".$}.

\vskip.1in
{\bf Definition:}  The sets $\mathcal B_r(f):= \mathcal B(\rho_r),$ $r=1.2\cdots \dim X,$ with the intervals $I$ converted into intervals with ends $\theta_i'$s and $(\theta_i +2\pi k)'$s and $\mathcal J_r(f):=\mathcal J(\rho_r)$  are the $r$-invariants of the map $f$. 
\vskip .1in
For a real valued map one has only
bar codes,  while for an angle valued map one has bar codes  and  Jordan cells.

We denote by $(V_r(f),T_r(f))$ the pair
$$
(V_r(f),T_r(f))=\bigoplus_{(\lambda,k)\in {\mathcal J}_r(f)}(\overline \kappa^k,T(\lambda,k))
$$ 
and refer to it as as the \emph{$r$-monodromy} of the angle valued map $f$.

\vskip .2in
\section {The main results}
\vskip .1in

Recall that for  $f:X\to \mathbb R$ a continuous map  denote by:  $X_t= f^{-1}(t)$ and  $X_{[t_1, t_2]}= f^{-1}[t_1,t_2].$  For $f:X\to S^1$ a continuous map denote by: 
$X_\theta= f^{-1}(\theta), \ \theta$ angle,   $\xi_f\in H^1(X;\mathbb Z)$ the cohomology class represented by $f$ and   $\tilde f: \tilde X\to \mathbb R$ the lift of $f$ to the infinite cyclic cover of $\tilde X\to X$   defined by $\xi_f.$ The covering $\tilde X\to X$ is the pull back of the  infinite cyclic cover $\mathbb R\to S^1$ by $f$
or any other map in the class $\xi_f.$
For $I$ an interval $\subset \mathbb R$ denote by $n_{\theta}(I)=\sharp\{k\in\mathbb Z\mid\theta+2\pi k\in I\}$ and for $J$ a Jordan cell  write $J=(\lambda(J), k(J)).$ 

\begin{theorem} \label{T1}\
 
\noindent1. If $f\colon X\to\mathbb R$ is a tame map then:

\begin{equation*}
\begin{aligned}
 &  \beta_r(X_t)\hskip .5in 
 &&=\sharp\bigl\{I\in\mathcal B_r (f)\bigm|I\ni t\bigr\}
\\ 
 &  \dim\rm {im} \bigl(H_r(X_t)\to H_r(X)\bigr)\hskip .5in
&&=\sharp\bigl\{I\in\mathcal B^c_r(f)\bigm| I \ni t\bigr\}\\
 &\beta_r(X)\hskip .5in
 &&=\sharp\mathcal B^c_r(f)+\sharp\mathcal B^o_{r-1}(f).
\end{aligned}
\end{equation*}
2. If $f\colon X\to S^1$ is a tame map then:

\begin{equation*}
\begin{aligned}
\hskip .5in 
 &\beta_r(X_\theta)
    &&=\sum_{I\in \mathcal B_r(f)}n_\theta(I)+\sum_{J\in\mathcal J_r(f)}k(J)
\\
\hskip .5in 
&\dim\rm {im}\bigl(H_r(X_\theta)\to H_r(X)\bigr)
&&=\sharp\bigl\{I\in \mathcal B^c_r(f)\bigm|\theta\in I\bigr\} + \sharp\bigl\{(\lambda,k)\in{\mathcal J}_{r}(f)\bigm|\lambda(J)=1\bigr\}\\
\hskip .5in 
&\beta_r(X)
&&=\begin{cases}
\sharp\mathcal B^c_r(f)+
\sharp\mathcal B^o_{r-1}(f)+\\
\sharp\bigl\{(\lambda,k)\in{\mathcal J}_{r}(f)\bigm|\lambda(J)=1\bigr\}+\\
\sharp\bigl\{(\lambda,k)\in{\mathcal J}_{r-1}(f)\bigm|\lambda(J)=1\bigr\}
\end{cases}\\
\hskip .5in 
&\beta N_r(X;\xi_f)&&=
\sharp\mathcal B^c_r(f)+\sharp\mathcal B^o_{r-1}(f).
\end{aligned}
\end{equation*}
\end{theorem}

\begin{theorem}\label {T2}    
 Let $f\colon X\to S^1$ be a tame map and  $\tilde {\mathcal B}(f): = \{ I'= I+2\pi k \mid k\in \mathbb Z, \  I\in \mathcal B(f)\}.$ Then:

1.
\begin{align*}
 & \beta_r(\tilde X_{[a,b]})
&&=\begin{cases} 
\sharp\{I'\in \tilde{\mathcal B}_r(f), I'\cap [a,b]  \rm {closed\  and} \ne\emptyset \}+\\
\sharp\{I'\in \tilde{\mathcal B}^o_{r-1}(f), I'\subset [a,b]\}+\\ 
\sum_{J\in \mathcal J_r(f)}k(J).
\end{cases}
\\
 & \dim\rm {im}\bigl(H_r(\tilde X_{[a,b]})\to H_r(\tilde X)\bigr) 
&&=\begin{cases} 
\sharp\{I \in\tilde{\mathcal B}^c_r(f), I\cap [a,b] \ne\emptyset\}+\\
\sharp\{I\in \tilde{\mathcal B}^o_{r-1}(f), I\subset [a,b]\}+\\ 
\sum_{J\in \mathcal J_r(f)} k(J).
\end{cases}
\\
 & \dim \rm{im}\bigl(H_r(\tilde X_{[a,b]}) \to H_r(X)\bigr) 
&&=\begin{cases}
\sharp\{I\in \mathcal B_r^c,[a,b]\cap (I+2\pi k)\ne 0\}+\\
\sharp\{I\in \mathcal B^o_{r-1}\mid   I+2\pi k \subset [a,b]\}+\\
\sharp\{J\in {\mathcal J}_{r}(f)| \lambda(J)=1\}+\\
\sharp\{J\in{\mathcal J}_{r-1}(f)| \lambda(J)=1\}. 
\end{cases}
\end{align*}

2.  
$V_r(\xi_f):=\ker (H_r(\tilde X)\to H^N_r(X;\xi_f))$ is a finite dimensional $\kappa$-vector space and $(V_r(\xi_f)\otimes \overline \kappa, T_r(\xi_f)\otimes \overline \kappa)=(V_r(f), T_r(f)).$ 

3. 
$H_r(\tilde X)=\kappa[t^{-1},t]^N\oplus V_r(\xi_f)$ as $\kappa[t^{-1}, t]$-modules with $N=\beta N_r(f)=\sharp{\mathcal B^c_r(f)} + \sharp{\mathcal B^o_{r-1}(f)}$.
\end{theorem}

Theorems \ref {T1} and \ref{T2} imply that for $f$ real valued  the number $\sharp \mathcal B_r^c+ \sharp \mathcal B^o_{r-1}$  is a  homotopy invariants and 
for $f$ angle valued the number $\sharp \mathcal B_r^c+ \sharp \mathcal B^o_{r-1}$ and the collection $\mathcal J_r(f)$ are homotopy invariants.
Therefore $C_r(f)$ can be regarded as points in the symmetric product $S^{\beta_r(X)}(R^2)$ resp. $S^{\beta N_r(X;\xi_f)}(\mathbb T)$  which are nice stratified spaces. Recall that $S^N(M)= ( M\times M\times \cdots M) / \Sigma_N$ where the product contains $N$ terms and $\Sigma_N$ denotes the $N-$symmetric group.

Let $C^0_{\rm{tame}} (M; \mathbb R)$ resp. $C^0_{\rm{tame}}(M; \mathbb S^1)$  denote the set of  tame maps with the topology induced from $C^0(M; \mathbb R)$ resp. $C^0(M; \mathbb S^1)$  equipped with the compact open topology. This set is dense in the space of all continuous maps.  

\begin{theorem}\label {T3}
The assignments $f  \rightsquigarrow C_r(f)$ is a continuous map on $C^0_{\rm{tame}}(M; \mathbb R)$ resp. $C^0_{\rm{tame}}(M; \mathbb S^1)$ hence has a continuous extension to the entire $C^0(M; \mathbb R)$ resp. $C^0_{\rm{tame}}(M; \mathbb S^1).$
\end{theorem}

As a consequence the configuration $C_r(f),$ hence the closed $r-$ bar codes and the open $(r-1)$ bar codes,  as well as the collection of Jordan cells can be defined for any continuous maps. Consequently  the monic polynomials
$P_r(f)(z)$  are well defined and the assignment $f \rightsquigarrow P_r(f)(z)$ continuous. Note that the collection $\mathcal J_r(f)$ remains constant on a connected component of $C^0(M; \mathbb S^1).$
Consequently, for $f:X\to \mathbb C\setminus 0$ a continuous map one has the monic polynomials $P_r(|f|)(z)$ and $Pr(f/ |f|)(z)$ which can be regarded as refinements of Betti numbers and Novikov-Betti numbers with respect to $f.$

The above results show that for $f:X\to S^1$ a tame map  only the bar codes in $\mathcal B^c_{\cdots}(f),$ $\mathcal B^o_{\cdots} (f)$ and the  Jordan cells  $J_{\cdots} (f)$ are relevant  for the topology of $X.$ 
The bar codes in $\mathcal B^{c,o}_{\cdots}$ and in $\mathcal B^{o,c}_{\cdots}$   are related only with the specifics of the map $f$ and have no contribution to the topology of $X.$ 
More about will be discussed in \cite {B12}.

\vskip .2in

\section {The  meaning of the bar codes}

\vskip .1in

For $f:X\to \mathbb R$ the following concepts are  fundamental to describe the  meaning of the invariants we have considered.
\begin{itemize}
\item The element $x\in H_r(X_t) $ is dead (to the right)  at $t'>t$ resp.  dead (to the left)  at $t"<t $ if its image by $H_r(X_t)\to H_r(X_{[t,t']})$ resp. by  $H_r(X_t)\to H_r(X_{[t",t]})$ vanishes.
\item The element $x\in H_r(X_t) $ is observable   at $t'\ne t$ r if its image by $H_r(X_t)\to H_r(X_{[t,t']})$ is contained in the image of $H_r(X_t')\to H_r(X_{[t,t']}).$
\end{itemize}

\begin{definition}\

\begin{enumerate}
\item For $x\in H_r(X_t)$ define $\tau^+(x)\in \mathbb R_+\cup \infty $  resp. $\tau^+(x)\in \mathbb R_+\cup \infty$ by the following property:
    $x$  is dead (to the right)at $t+\tau^+(x)$ resp. (to the left) at $t-\tau^-$ but not before, i.e for $t'$ with $t<t'< t+\tau^+(x)$ resp.  $t>t' >t-\tau^-.$
\item For $x\in H_r(X_t)$define $o^+(x)\in \mathbb R_+\cup \infty $ resp $o^-(x)\in \mathbb R_+\cup \infty $ \footnote{if $X$ is compact in particular if $f$ is tame as defined above $o^\pm(x)\in \mathbb R$, so can not be infinite } by the following property:
$x$ is  observable  at  $t+o^+(x)$ resp. $t-o^-(x)$ but not at $t+o^+(x)+\epsilon$  resp. $t- o^-(x)-\epsilon$ for $\epsilon >0.$
\end{enumerate} 
\end{definition}
 
 \vskip .2in

\begin{definition}
For $f:X\to \mathbb R,$  with critical values $s_i$ 
denote by: 
\begin{enumerate}
\item  $N_r(s_i, s_j)$ ($s_i\leq s_j$) the maximal number of linearly independent elements $x\in H_r(X_t)$ with $t+\tau^+(x)= s_j, t-\tau^-(x)= s_i$ for  any  $t$  in the open interval $(s_i, s_j),$
\item  $N_r[s_i, s_j]$ ($s_i< s_j$) the maximal number of linearly independent elements $x\in H_r(X_t)$ with $t+o^+(x)= s_j, t-o^-(x)= s_i$ for   any $t$ 
\footnote{ it suffices to happen for one $t$ and then it happens for any other $t$} in the open interval $(s_i, s_j),$
\item   $N_r(s_i, s_j]$ ($s_i< s_j$)the maximal number of linearly independent elements $x\in H_r(X_t)$ with $t+o^+(x)= s_j, t-tau^-(x)= s_i$ for  any $t$ in the open interval $(s_i, s_j),$
\item  $N_r[s_i, s_j)$ For $s_i< s_j$ the maximal number of linearly independent elements $x\in H_r(X_t)$ with $t+\tau^+(x)= s_j, t-o^-(x)= s_i$ for  any $t$ in the open interval $(s_i, s_j).$
\end{enumerate}
\end{definition}
For $f$ real values the number $N_r\{s_i, s_j\}$ represents the multiplicity of the bar code $\{s_i, s_j\}$ with the convention that  non existence of such bar codes means multiplicity zero.

For $f$ angle valued with critical values $\theta_i$ the number $N_r\{\theta_i, \theta_j +2\pi k\}$ represents  the multiplicity of the bar code $\{\theta_i, \theta_j+ 2\pi k \}$  which is the same as the multiplicity of  $\{\theta_i, \theta_j+ 2\pi k \}$ for the real valued $\tilde f :\tilde X\to \mathbb R$ of $f$.

Note that for $f:X\to S^1$ and $\theta\in (0, 2\pi]$ one can have elements $x\in H_r(X_\theta)=H_r(\tilde X_{\theta})$ which never die and remain observable for ever.  The existence of such elements is guarantied by the presence of Jordan cells . The Jordan cells provide rather complete information on the maximal number such elements which remain observable and linearly independent for any $\theta'$  as well as about how they return in $H_r(\tilde X_{\theta+2\pi})= H_r(\tilde X_{\theta})$ when $\theta'$    goes from $\theta$ to $\theta +2\pi,$ equivalently how do they change when observed in $H_r(X_{\theta+2\pi}).$ 
\vskip .2in

\section {About the proof (the canonical long exact sequence)}
\vskip .1in
Since a tame real valued map can be regarded as a tame angle valued map (by identifying $\mathbb R$ to  an open subset of 
$S^1,$)  we will consider only  the case of tame angle valued maps.

Let $f\colon X\to S^1$ be a tame map with $m$ critical angles $\theta_1,\theta_2,\dotsc,\theta_m$ and regular angles $t_1,t_2,\dotsc,t_m$.
First observe that, up to homotopy, the space $X$ and the map $f\colon X\to S^1$ can be 
regarded as the iterated mapping torus $\mathcal T$ and the map  $f^{\mathcal T}: \mathcal T \to[0,m]/{\sim}$ described below. 
Consider the collection of spaces and continuous maps:

{\scriptsize
\begin{center}
$$
\xymatrix{ &X_{\theta_1}\\
X_{t_1}\ar[ur]_{b_1}\ar[d]^{a_2}& & X_{t_m}\ar[ul]^{a_1}\ar[d]_{b_m}\\
X_{\theta_2}& & X_{\theta_m}\\
X_{t_{m-2}} \ar@{<.>}[u]\ar[dr]^{a_{m-1}}& & X_{t_{m-1}}\ar[u]^{a_m}\ar[dl]_{b_{m-1}}\\
&X_{\theta_{m-1}}
}
$$
\end{center}
}
\noindent with $R_i:=X_{t_i}$ and $X_i:=X_{\theta_i}.$ Denote by $\mathcal T=T(\alpha_1 \cdots \alpha_m;\beta_1\cdots\beta_m)$ the space obtained from the disjoint union 
$$
\Bigl(\bigsqcup_{1\leq i\leq m}R_i\times[0,1]\Bigr)\sqcup\Bigl(\bigsqcup_{1\leq i\leq m}X_i\Bigr)
$$
by identifying $R_i\times\{1\}$ to $X_i$ by $\alpha_i$ and $R_i\times\{0\}$ to $X_{i-1}$ by $\beta_{i-1}$. 
Denote by $f^{\mathcal T}\colon\mathcal T\to[0,m]/{\sim}=S^1$ the map given by$f^{\mathcal T}\colon R_i\times[0,1]\to[i-1,i]$ is the projection on
$[0,1]$ followed by the translation of $[0,1]$ to $[i-1,i]$ and $ [0,m]/{\sim}$ the space obtained from the segment $[0,m]$ by identifying the ends. 
The map $f^{\mathcal T}\colon\mathcal T\to[0,m]/{\sim}$ is a {\it homotopical reconstruction} of $f\colon X\to S^1$ provided that, with the choice of angles  
$t_i$, $\theta_i,$ the maps $a_i$, $b_i$ are those described in section~\ref{S6} for $X_i:=f^{-1}(\theta_i)$ and  $R_i:=f^{-1}(t_i)$.

Let $\mathcal P'$ denote the space obtained from the disjoint union 
$$
\Bigl(\bigsqcup_{1\leq i\leq m}R_i\times(\epsilon,1]\Bigr)\sqcup\Bigl(\bigsqcup_{1\leq i\leq m} X_i\Bigr)
$$
by identifying $R_i\times\{1\}$ to $X_i$ by $\alpha_i$, and $\mathcal P''$ denote the space obtained from the disjoint union 
$$
\Bigl(\bigsqcup_{1\leq i\leq m}R_i\times[0,1-\epsilon\Bigr)\sqcup\Bigl(\bigsqcup_{1\leq i\leq m}X_i\Bigr)
$$
by identifying $R_i\times\{0\}$ to $X_{i-1}$ by $\beta_{i-1}$.

Let $\mathcal R=\bigsqcup_{1\leq i\leq m}R_i$ and $\mathcal X=\bigsqcup_{1\leq i\leq m}X_i$. Then, one has:
\begin{enumerate}
\item 
$\mathcal T=\mathcal P' \cup \mathcal P''$,
\item 
$\mathcal P'\cap\mathcal P''=\bigl(\bigsqcup_{1\leq i\leq m}R_i\times(\epsilon,1-\epsilon)\bigr)\sqcup\mathcal X$, and
\item 
the inclusions $\bigl(\bigsqcup_{1\leq i\leq m}R_i\times\{1/2\}\bigr)\sqcup\mathcal X\subset\mathcal P'\cap\mathcal P''$
as well as the obvious inclusions $\mathcal X\subset\mathcal P'$ and $\mathcal X\subset\mathcal P''$ are homotopy equivalences.
\end{enumerate}
The Mayer--Vietoris long exact sequence applied to $\mathcal T= \mathcal P'\cup \mathcal P''$ leads to the diagram:
\begin{center}
$$
\xymatrix{
&                                                       & H_r(\mathcal R)\ar[r]^{M (\rho_r)}                                                  & H_r(\mathcal X)\ar[rd]                                               &&\\
\cdots\ar[r]&H_{r+1}(\mathcal T)\ar[ur]
\ar[r]^-{\partial_{r+1}}&H_r(\mathcal R)\oplus H_r(\mathcal X)\ar[u]^{pr_1}\ar[r]^N  &H_r(\mathcal X)\oplus H_r(\mathcal X)\ar[u]^{(Id,-Id)}{\ar[r]^-{(i^r, -i^r)}} &H_r(\mathcal T)\ar[r] &\\
&                                                       &H_r(\mathcal X)\ar[u]^{in_2}\ar[r]^{Id}                                           &H_r(\mathcal X)\ar[u]^{\Delta}                                                 &&
}
.$$
Diagram~2
\end{center}

Here $\Delta$ denotes the diagonal, $in_2$ the inclusion on the second component, $pr_1$ the projection on the first component, 
$i^r$ the linear map induced in homology by the inclusion $\mathcal X\subset \mathcal T.$ The matrix  $M(\rho_r)$ is defined by 
\begin{equation*}
M(\rho_r)= \begin{pmatrix}
\alpha^r_1 & -\beta^r_1 & 0          & \cdots         & 0\\
0          & \alpha^r_2 & -\beta^r_2 & \ddots         & \vdots\\
\vdots     & \ddots     & \ddots     & \ddots         & 0 \\
0          & \cdots     & 0          & \alpha^r_{m-1} & -\beta^r_{m-1}\\
-\beta^r_m & 0          & \cdots     & 0              & \alpha^r_m
\end{pmatrix}_. 
\end{equation*}
with $\alpha^r_i\colon H_r(R_i)\to H_r(X_i)$ and $\beta^r_i\colon H_r(R_{i+1})\to H_r(X_i)$ induced by the maps $\alpha_i$ and $\beta_i$ and the matrix $N$ is defined by

\begin{equation*}
\begin{pmatrix}
\alpha^r & Id \\
-\beta^r &  Id
\end{pmatrix}
\end{equation*}


where $\alpha^r$ and $\beta^r$ are the matrices 
$$
\begin{pmatrix}
\alpha^r_1 & 0          & \cdots &    0 \\
0          & \alpha^r_2 & \ddots &\vdots \\
\vdots     & \ddots     & \ddots &0 \\
0          & \cdots     & 0      & \alpha^r_{m-1} 
\end{pmatrix}
\quad\text{and}\quad
\begin{pmatrix}
0         & \beta^r_1 & 0         & \dots  & 0 \\
0         & 0         & \beta^r_2 & \ddots & \vdots \\
\vdots    & \vdots    & \ddots    & \ddots & 0 \\
0         & 0         & \dots     & 0      & \beta^r_{m-1}\\
\beta^r_m & 0         & \dots     & 0      & 0
\end{pmatrix}_.
$$

The long exact sequence 
\begin{equation}\label{MV}
\boxed{\cdots\to H_r(\mathcal R)\xrightarrow{M(\rho_r)}H_r(\mathcal X)\to H_r(\mathcal T)\to H_{r-1}(\mathcal R)\xrightarrow{M(\rho_{r-1})}H_{r-1}(\mathcal X)\to\cdots}
\end{equation}
derived from  Diagram~2 is referred to as the {\it canonical sequence} associated with a tame. 

This long exact sequence implies the short exact sequence

%
\begin{equation}\label{EE13}
0\to\rm{coker} M(\rho_r)\to H_r(\mathcal T)\to\ker M(\rho_{r-1})\to0
\end{equation}
and then the  noncanonical isomorphism
\begin{equation} \label{E11}
H_r(\mathcal T)=\rm{coker} M(\rho_r)\oplus\ker M(\rho_{r-1}).
\end{equation}
Any splitting $s\colon\ker M(\rho_{r-1})\to H_r(\mathcal T)$ in the short exact sequence~\eqref{EE13} provides an isomorphism~\eqref{E11}.
The calculation of $\ker M(\rho)$ and $\rm{coker} M(\rho)$ for $\rho= \rho_r$ is reduced to the case $\rho$ is indecomposable  hence to bar codes and Jordan cells cf \cite {BD11} calculation provided in \cite{BD11} Proposition 5.3. 

Note that the long exact sequence \eqref{MV} holds also for homology with local coefficients (i.e. homology with coefficients in a representation of the fundamental group). 
Such sequence   can be derived from a similar diagram as  Diagram 2, where instead of homology with coefficients in $\kappa$ one uses homology with local coefficients.
Of particular interest is the case the local coefficients system is   
$u\xi_f,$ the representation described by the composition $H_1(M;\mathbb  Z) \xrightarrow{\xi_f}\mathbb Z \xrightarrow {\hat u}\mathbb C^\ast= \mathbb C\setminus 0$ with $\hat u(n)= u^n,\ u\in \mathbb C\setminus 0. $  
 In this case the vector spaces $H_r(\mathcal R)$ and $H_r(\mathcal X)$ are independent on $u\xi_f$ and represent always the cohomology with coefficients in the trivial representation (corresponding to $u=1)$ hence with coefficients in the field $\kappa.$ 
Manipulations of these sequences  cf \cite{BD11}, \cite {BH12} lead to the proof of Theorems \ref{T1}, \ref{T2}. The proof of Theorem \ref{T3} requires refinements of Theorems \ref{T1} and \ref{T2} and will be contained in \cite{B12}.

Note:

1.The canonical long exact sequence contains more information than used in the present discussion. In case $\kappa$ is a field of characteristic zero $H_r(\mathcal R)$ and $H_r(\mathcal X)$contain inside the lattice of integral homology. This can be used for calculating a more subtle invariant  "torsion". 

2.Theorems \ref{T1}, \ref{T2} imply that up to  an isomorphism (of vector spaces) the canonical long exact sequence it is completely determined by the bar codes and the Jordan cells.
\vskip .2in

\section {About computability of the bar codes and the Jordan cells}
\vskip .1in
For $f:X\to \mathbb R$ or $f:X\to S^1$  simplicial  maps  algorithms of relative low complexity   to calculate the bar codes and Jordan cells are described in   \cite{BD11}.
Here {\it simplicial}  means that $X$ is a simplicial complex and, when the target of $f$ is $\mathbb R,$  the  restriction of $f$ to any simplex  $\sigma$  of $X$ is linear, and when the target is $S^1,$ any lift  $\tilde f:\sigma\to \mathbb R$ of $f |_\sigma \to S^1$ ($\pi \cdot \tilde f= f|_\sigma$ with $\pi:\mathbb R\to S^1$ the universal cover) is linear.
Note that any simplicial map $f$ is tame and its critical values are among the values of $f$ on the vertices of $X.$ 

The algorithms we proposed consist of two steps.  In the first step one inputs the simplicial complex and the values of $f$ on the vertices of $X$ and derive which of these values are critical and then by choosing regular values (for example midle between two consecutive critical values)  the representations  $\rho_r$ as  collections of matrices $\{\alpha_i^r, \beta_i^r\}.$ 
A summary presentation of this part is provided in \cite{BD11}. The second part inputs the representations $\rho_r$ and outputs  the  bar codes and the Jordan cells.
Details  are provided in  the Appendix to \cite {BD11}. 

As long as the first input is concerned, to record the simplicial complex $X$ we choose a total order of the set of vertices, $\{v_1, v_2, \cdots \},$ extend this  order to a total order of the set of all simplices of $X$
such that the following two conditions hold:
 
 C1: If $\tau$ is a face of $\sigma$ then $\tau< \sigma,$
 
 C2: If $\dim \tau  <\dim \sigma$ then $\tau<\sigma.$ It is obvious that such total orders exists.  

Note that the ordering of the vertices provide an orientation on each simplex. 

One records the simplicial complex $X$  as an $N \times N$ upper triangular matrix $M(X)$ with zero on diagonal and entries $0, +1, -1$ where $N$ the cardinality of the set of all simplices. Precisely the entry corresponding to the pair $(\tau, \sigma)$ is zero if $\tau$ is not a face of $ \sigma $ and equal to $\pm 1$ if it is. In this case   is $+1$ if the orientation of $\sigma$ induces  the orientation of $\tau$ and $-1$ otherwise.
 
 To the matrix called $M(X)$ one add the values of $f$ on vertices. In the first step one determines which values of $f$ are critical  and then using sub matrices of  $M(X)$ (possibly enhanced)  one recover the matrices $M(\rho_r)$ equivalently  the representations $\rho_r$  as indicated in \cite{BD11}.

In the second step a new algorithm whose input is  the matrix $M(\rho_r)$ describing the representations $\rho_r$ and output is the barcodes and the Jordan cells  finalize the calculations.  
More details about this algorithm can be found in \cite{BD11}.
\vskip .2in

\section {Examples}
\vskip .1in

The Picture below 
describes a tame real valued map $p: Y\to [0,2\pi]\subset \mathbb R$ and an angle valued map $f\colon X\to \mathbb S^1$ whose bar codes and Jordan cells are given in the  tables below. 
\begin{figure}[h]
\includegraphics{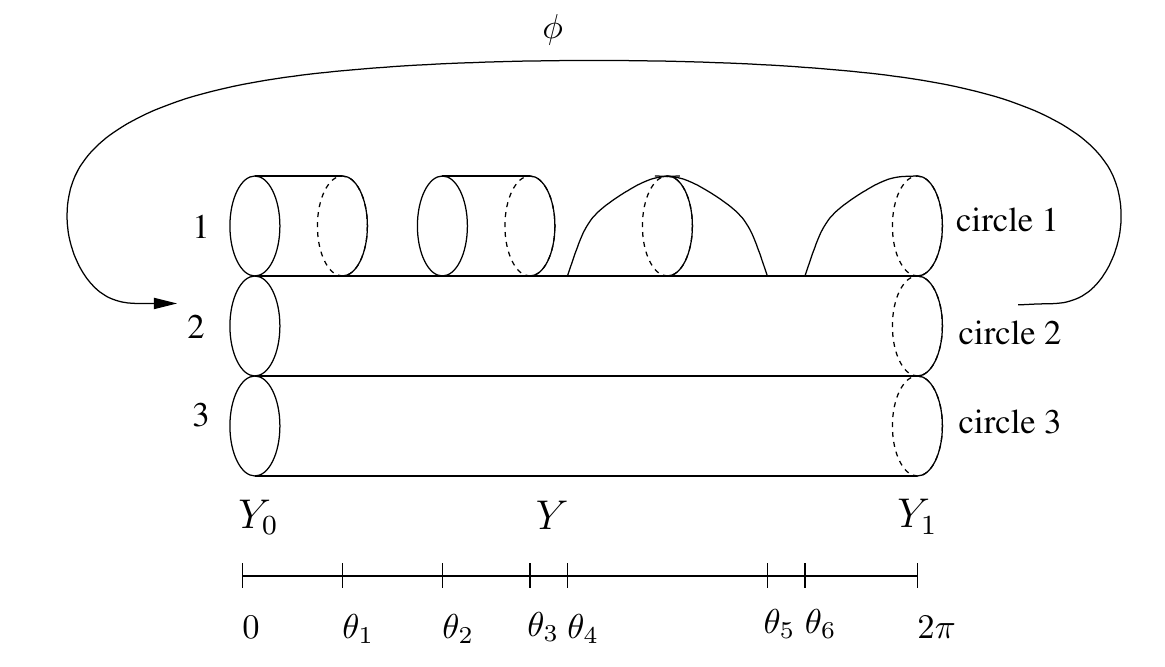}
\caption{The tames maps $p$ and $f$.}

\label{cv-ex}
\end{figure}
The space $X$ is obtained from $Y$  by identifying its  
right end $Y_1$ (a union of three circles) to the left end 
$Y_0$ (a union of three circles) following 
the map $\phi\colon Y_1\to Y_0$  described as follows:

- circle 1 goes 3 times around circle 1

- circle 2 go 2 time around circle 2 

- circle 3 goes1 time around 2 counter clockwise and 2 times around circle 3.

\noindent The map 
$p:Y\to \mathbb R$ is the projection on $[0,2\pi]$ and the map $f\colon X\to S^1$ is induced by 
the projection of $p: Y \to [0,2\pi]$  by passing to the quotient spaces $X=Y/\varphi$ and $[0,2\pi]/\sim.$ 
Note that  $H_1(Y_1)=H_1(Y_0)=\kappa\oplus\kappa\oplus\kappa$
and $\phi$ induces a linear map in $H_1$-homology represented by the matrix 
\begin{equation*}
\begin{pmatrix}
3&0 & 0\\
1&2&-1\\
0 & 0&2            
\end{pmatrix}.
\end{equation*}

The bar codes of the map $p$ are given in the Table 1. There are no bar codes in dimension 2 since each fiber of $f$ is one-dimensional.
\vskip .1in

\hskip 1.5in
{\scriptsize
\begin{tabular}{c|c} 
& $r$-invariants for $p$\\ \hline
\begin{tabular}{l}
\end{tabular} & 
\begin{tabular}{c|c|c}
dimension & bar codes  \\ \hline
0  &  $[0,2\pi] $\\ \hline
   & $[0,  \theta_1]$\\
   &$(\theta_3, \theta_5)$\\
   &$(\theta_6, 2\pi]$\\
1 & $[\theta_2,\theta_3]$ & \\
   & $[0, 2\pi] $  &  \\
   & $[0, 2\pi] $  &  \\
\end{tabular} \\ \hline 
\end{tabular}}
\vskip .1in
\hskip 2 in { Table 1: } 

\vskip .1in

For the angle valued map $f:X\to S^1$ 
there are no bar codes or Jordan cells in dimension 2 since each fiber of $f$ is one-dimensional 
and, as all fibers are connected in dimension zero we have only one Jordan cell $\rho^{II}(1; 1).$
The bar codes and the Jordan cells in dimensions $0$ and  $1$ are described Table 2.in  More details on their calculation are 
presented in  \cite {BD11} and \cite {BH12}.

\vskip .1in

\hskip 1in
{\scriptsize
\begin{tabular}{c|c}
$r$-invariants of $f$ & \\ \hline
\begin{tabular}{l}
\end{tabular} 
\begin{tabular}{c|c|c}
dimension & bar codes & Jordan cells \\ \hline
0  &  & $(1,1)$ \\ \hline
    & $(\theta_6, \theta_1+2\pi]$ & $(2,2)$ \\
1  & $[\theta_2,\theta_3]$ & \\
    & $(\theta_4,\theta_5)$&  \\
\end{tabular} \\ \hline 
\end{tabular}}
\vskip .1in
 
 \hskip 2in
 Table 2.

\vskip.2in






\end{document}